\documentclass[journal ]{new-aiaa}
\usepackage[utf8]{inputenc}
\usepackage{textcomp}

\usepackage{graphicx}
\usepackage{amsmath}
\usepackage{mathtools}
\usepackage{siunitx}
\usepackage{longtable,tabularx}

\newcommand{\br}{\textbf{\emph{r}}}

\newcommand{\bu}{\textbf{\emph{u}}}
\newcommand{\bv}{\textbf{\emph{v}}}
\newcommand{\be}{\textbf{\emph{e}}}

\newcommand{\bs}{\textbf{\emph{s}}}

\newcommand{\bp}{\textbf{\emph{p}}}
\newcommand{\bq}{\textbf{\emph{q}}}
\newcommand{\bw}{\textbf{\emph{w}}}
\newcommand{\bx}{\textbf{\emph{x}}}

\newcommand{\bC}{\textbf{\emph{C}}}
\newcommand{\bD}{\textbf{\emph{D}}}
\newcommand{\bN}{\textbf{\emph{N}}}
\newcommand{\bS}{\textbf{\emph{S}}}
\newcommand{\bT}{\textbf{\emph{T}}}

\newcommand{\bell}{\boldsymbol{\ell}}

\newcommand{\bxi}{\boldsymbol{\xi}}

\newcommand{\PP}{\mathbb{P}}

\title{Solving the Gibbs Problem with Algebraic Projective Geometry}

\author{Michela Mancini \footnote{Graduate~Research~Assistant,~Guggenheim~School~of~Aerospace Engineering,~Atlanta,~30332,~GA} and John A. Christian \footnote{Associate Professor, ,~Guggenheim~School~of~Aerospace Engineering,~Atlanta,~30332,~GA,~Associate~Fellow~AIAA}}
\affil{Georgia~Institute~of Technology,~Atlanta,~30332,~GA.}

\begin{document}

\maketitle

\section{Introduction}
Orbit determination (OD) from three position vectors is one of the classical problems in astrodynamics. Early contributions to this problem were made by J. Willard Gibbs in the late 1800s \cite{Gibbs:1889} and OD of this type is known today as ``Gibbs Problem.'' There are a variety of popular solutions to the Gibbs problem. While some authors solve for the orbital elements directly \cite{Gurfil:2016}, most contemporary discussions are based on a vector analysis approach inspired by Gibbs himself. An especially nice version of this was provided by Bate, Mueller, and White \cite{BMW} and has since become a commonly adopted presentation in astrodynamics textbooks \cite{Curtis:2020,Vallado:2007,Spencer:2023}.

Omitting the derivation (see any of the textbooks just mentioned \cite{BMW,Curtis:2020,Vallado:2007,Spencer:2023}), the usual solution to the Gibbs problem taught to students of aerospace engineering is approximately as follows. Consider a body in a Keplerian orbit. Suppose that the analyst only has access to three distinct position vectors $\{\br_i\}_{i=1}^3$ describing the location of the orbiting body at three different times. Solving for the orbit requires determination of the unknown velocities $\{\bv_i\}_{i=1}^3$. This task may be accomplished by first computing the three intermediate vectors $\bN$, $\bS$, and $\bD$
\begin{subequations}
\begin{equation}
    \bN =   r_1 (\br_2 \times \br_3)
          + r_2 (\br_3 \times \br_1)
          + r_3 (\br_1 \times \br_2)
\end{equation}
\begin{equation}
    \bD = (\br_1 \times \br_2) + (\br_2 \times \br_3) + (\br_3 \times \br_1)
\end{equation}
\begin{equation}
    \bS =   ( r_2 - r_3 ) \br_1
          + ( r_3 - r_1 ) \br_2
          + ( r_1 - r_2 ) \br_3
\end{equation}
\end{subequations}
where $r_i = \| \br_i \|$ are the known distances from the central body. Assuming the central body has a gravitational constant of $\mu$, the unknown velocities $\bv_i$ may be directly computed as
\begin{equation}
    \bv_i = \sqrt{ \frac{\mu}{ND} } \left( \frac{\bN \times \br_i}{r_i} + \bS \right)
\end{equation}
where $N = \| \bN \|$ and $D = \| \bD \|$. The Gibbs problem has now been solved.

An alternate, but equivalent, way to think about this problem is to use the vectors $\bN$, $\bS$, and $\bD$ to compute the orthonormal basis vectors of the perifocal frame $(\hat{\bp},\hat{\bq},\hat{\bw})$ that describe the orientation of the orbit. These basis vectors are computed as
\begin{subequations}
\begin{equation}
    \hat{\bp} = (\bS \times \bN) / (S N)
\end{equation} 
\begin{equation}
    \hat{\bq} = \bS / S
\end{equation} 
\begin{equation}
    \hat{\bw} = \bN / N
\end{equation} 
\end{subequations}
The ordered bases $(\hat{\bp},\hat{\bq},\hat{\bw})$ clearly form a right-handed system. Now, with the perifocal frame defined, one may also find the semi-latus rectum $p$ and eccentricity $e$ that describe the size and shape, respectively, of the orbit
\begin{subequations}
\begin{equation}
    p = N/D
\end{equation}
\begin{equation}
    e = S/D
\end{equation}
\end{subequations}
The Gibbs problem has once again been solved.

This work presents a completely different solution to those just described. Although there is nothing wrong with the vector analysis approach, some interesting insights may be gained  by considering the problem from the perspective of algebraic projective geometry. Such an algebraic solution is presented here. The OD procedure is based upon a novel and geometrically meaningful solution to the algebraic fitting of an ellipse with a focus at the origin using only three points. Although the final OD result is identical to the classical vector analysis approach pioneered by Gibbs, this new algebraic solution is interesting in its own right.

\section{An Algebraic View of Conics}\label{sec:an_algebraic_view}
Keplerian orbits are planar and follow a path described by a conic (i.e., a circle, ellipse, parabola, or hyperbola), with the gravitating body at the one of the conic's foci. Thus, within the orbital plane, the orbit is simply a two-dimensional conic.

As was shown by Descartes, the algebraic description of a conic is the locus of points $(x,y)$ satisfying a polynomial of degree two in two variables
\begin{equation}
    \label{eq:ConicPoly}
    A x^2 + B x y + C y^2 + Dx + Ey + G = 0
\end{equation}
If one chooses to work in two-dimensional projective space $\PP^2$, then points (not at infinity) on the projective plane may be written as $\bar{\bx} = [x;y;1]$. Thus, the conic locus is given by
\begin{equation}
    \label{eq:ConicLocusConstraint}
    \bar{\bx}^T \bC \bar{\bx} = 0
\end{equation}
where $\bC$ is the $3 \times 3$ symmetric matrix of ambiguous scale
\begin{equation}
    \label{eq:ConicLocusMatrix}
    \bC \propto
    \begin{bmatrix}
          A & B/2 & D/2 \\
          B/2 & C & E/2 \\
          D/2 & E/2 & F
    \end{bmatrix}
\end{equation}
Now, the equation for a line is
\begin{equation}
    a x + b y + c = 0
\end{equation}
Thus, defining a line in $\PP^2$ as $\bell = [a;b;c]$, a point $\bar{\bx}$ lies on the line $\bell$ if and only if $\bar{\bx}^T \bell = 0$. It follows directly from this definition and Eq.~\eqref{eq:ConicLocusConstraint} that, if $\bx$ is a point on the conic, then the line $\bell = \bC \bar{\bx}$ must be tangent to the conic since
\begin{equation}
    \bar{\bx}^T \bell = \bar{\bx}^T \bC \bar{\bx} = 0
\end{equation}
Moreover, given a tangent to the conic of $\bell = \bC \bar{\bx}$, one may also write $\bar{\bx} = \bC^{-1} \bell$. This is always possible for a non-degenerate orbit since $\bC$ is always full rank for a non-degenerate conic. Substitution of this result into Eq.~\eqref{eq:ConicLocusConstraint} yields
\begin{equation}
    \bar{\bell}^T \bC^{-1} \bar{\bell} = 0
\end{equation}
which describes all the lines tangent to the conic. This is referred to as the conic envelope.

Keplerian orbits are not arbitrary conics. Instead, Keplerian orbits are specifically those conics having a focus at the gravitating body. Without loss of generality, the present analysis will be made easier by placing the gravitating body at the origin such that the conic has a focus at the origin. How to write this constraint directly in terms of the coefficients $A,B,C,D,E,F$ from Eq.~\eqref{eq:ConicPoly} and Eq.~\eqref{eq:ConicLocusMatrix} is not always straightforward to see. However, a conic locus with a foci at the origin must have a conic envelope that is a circle \cite{Chakerian:2001}. Recognizing that the velocity vector is always tangent to the conic, this statement turns out to be equivalent to the orbital hodograph being a circle---a fact first demonstrated by W. R. Hamilton \cite{Hamilton:1847}.

If the conic envelope is a circle, then the matrix $\bC^{-1}$ has a simple form of
\begin{equation}\label{eq:Cinvmatrix}
    \bC^{-1} \propto \begin{bmatrix}
    1 & 0 & X\\
    0 & 1 & Y\\
    X & Y & -Z^2
    \end{bmatrix}
\end{equation}
which may be inverted to find a simple expression for the conic locus constrained to have one of its foci at the origin
\begin{equation}\label{eq:Cmatrix}
    \bC\propto \begin{bmatrix}
        -(Y^2+Z^2) & XY & -X\\
        XY & -(X^2+Z^2) & -Y\\
        -X & -Y & 1
    \end{bmatrix}
\end{equation}

Now, using Eqs.~\eqref{eq:ConicLocusMatrix}, Eq.~\eqref{eq:Cinvmatrix}, and Eq.~\eqref{eq:Cmatrix} in conjunction with the identities from Ref.~\cite{Mancini:2024}, the geometric parameters of the orbit may be written in terms of the six polynomial coefficients $A,\hdots,F$ or in terms of the three parameters $X,Y,Z$.

The direction from the focus to the periapsis is given by the unit vector 
\begin{equation}
    \hat{\bp} = \frac{1}{\sqrt{D^2 + E^2}}\begin{bmatrix}
        -D \\ -E \\0
    \end{bmatrix}
    =\frac{1}{\sqrt{X^2 + Y^2}}\begin{bmatrix}
        X \\ Y \\ 0
    \end{bmatrix}
\end{equation}
where, of course, $\hat{\bp}$ is expressed in the same frame as the measured 2-D positions (x,y). Thus, it becomes clear that the parameters $X,Y$ in Eq.~\eqref{eq:Cinvmatrix} and Eq.~\eqref{eq:Cmatrix} describe the direction from the focus (at the origin) to periapsis. This fact gives helpful geometric meaning to derivations that follow.

The semi-major axis $a$ and semi-minor axis $b$ may be computed as
\begin{equation}
    a = -2\frac{\sqrt{D^2F^2-4AF^3}}{E^2+4AF}=\frac{\sqrt{X^2+Y^2+Z^2}}{Z^2}
\end{equation}
\begin{equation}
    b = 2\sqrt{-\frac{F^2}{E^2+4AF}} =\frac{1}{Z}
\end{equation}
and the eccentricity $e$ may be computed as
\begin{equation}
    \label{eq:ecc}
    e = \sqrt{\frac{E^2+D^2}{D^2-4AF}}= \sqrt{\frac{X^2+Y^2}{X^2+Y^2+Z^2}}
\end{equation}
Finally, the semi-latus rectum $p$ may be computed as
\begin{equation}
    \label{eq:SemiLatusRectum}
    p =2\sqrt{\frac{F^2}{D^2-4AF}}= \frac{1}{\sqrt{X^2+Y^2+Z^2}}
\end{equation}

\section{ Orbit Determination}
It is possible to use the conic relations from Section~\ref{sec:an_algebraic_view} to algebraically solve the Gibbs problem. The procedure has two parts. The first part involves finding the orbit plane, and the second part involves fitting a conic to the two-dimensional (2-D) positions within the orbit plane.

\subsection{Finding the Orbit Plane}
Keplerian orbits are planar, and it follows that the three position vectors $\{\br_i\}_{i=1}^3$ are coplanar. The normal to the orbit plane $\bw$ may be found as the cross-product of any two of the position vectors
\begin{equation}
    \bw \propto \br_1 \times \br_2 \propto \br_2 \times \br_3 \propto \br_3 \times \br_1
\end{equation}
or as the solution to a null space problem
\begin{equation}
    \begin{bmatrix}
        \br_1^T \\ \br_2^T \\ \br^T_3
    \end{bmatrix} \bw = \textbf{0}_{3 \times 1}
\end{equation}
No matter the method for finding $\bw$, the next task is to define two basis vectors $\be_1$ and $\be_2$ that span the orbit plane. These may be used to form a temporary frame to perform OD calculations. Since the orientation of these basis vectors is arbitrary, any principled choice will suffice. In this work, the two basis vectors are computed as
\begin{equation}
   \be_2 = \frac{\bw \times \br_1}{ \| \bw \times \br_1 \|}
\end{equation}
and
\begin{equation}
    \be_1 = \be_2 \times \bw
\end{equation}
This may be used to form the right-handed coordinate frame $N$ with axes $(\be_1,\be_2,\bw)$. The attitude transformation matrix from the inertial frame to this new intermediate frame $\bT^I_N$ may be computed as
\begin{equation}
    \bT^I_N = 
    \begin{bmatrix}
        \be_1^T \\ \be_2^T \\ \bw^T
    \end{bmatrix}
\end{equation}
Thus, the position vectors rotated into the intermediate frame $N$ are zero in the third coordinate axis
\begin{equation}
    \begin{bmatrix}
        x_i \\ y_i \\ 0
    \end{bmatrix} = 
    \bT^I_N \br_i 
\end{equation}
where is evident that $r_i = \| \br_i \| = \sqrt{x^2_i + y^2_i}$. The homogeneous coordinate ($\bar{\bx} \in \PP^2$) of each 2-D position vector may then be written as
\begin{equation}
    \label{eq:PointsToOrbitPlane}
    \bar{\bx}_i \propto
    \begin{bmatrix}
        x_i \\ y_i \\ 1
    \end{bmatrix} = 
    \begin{bmatrix}
        1 & 0 & 0 & 0 \\
        0 & 1 & 0 & 0 \\
        0 & 0 & 0 & 1 \\
    \end{bmatrix}
    \begin{bmatrix}
        \bT^I_N \br_i \\ 1
    \end{bmatrix}
\end{equation}

\subsection{Algebraic Conic Fitting}
Conic fitting generally requires five points. As a consequence, most algebraic (and numeric) conic fitting algorithms require at least five points \cite{Fitzgibbon:1999,Kanatani:2006,AlSharadqah:2012,Szpak:2012}. None of these conventional methods are applicable to the algebraic Gibbs problem.

Rather than fitting a general conic to five points, the algebraic Gibbs problem requires the fit of a conic with an origin at the focus to three points. The constraint that one of the foci is at the origin removes two degrees of freedom from the unconstrained five-point ellipse fitting problem. The authors are unaware of any existing published solutions to three-point ellipse fitting problem with a constrained focus location. Thus, a method for solving this problem is developed here.

Equation~\eqref{eq:PointsToOrbitPlane} may be used to transform the original position vectors $\{\br_i\}_{i=1}^3$ into points in $\PP^2$ within the orbit plane. Once in the orbit plane, the relation between the conic $\bC$ and the points $\{\bar{\bx}_i \}_{i=1}^3$ is given by Eq.~\eqref{eq:ConicLocusConstraint}. If $\bC$ has the specific structure of Eq.~\eqref{eq:Cmatrix}, then the focus is at the origin. Recognizing that the resulting equation is linear in $Z^2$, this may be solved for $Z^2$ using any of the points $\{\bar{\bx}_i \}_{i=1}^3$. For example, using the first point $\bar{\bx}_1$, the value of $Z^2$ may be computed directly as
\begin{equation}\label{eq:Z}
    Z^2=\left(\frac{1}{r_1^2}\right)\bar{\bx}_1^T\begin{bmatrix}
        -Y^2 & XY & -X\\
        XY & -X^2 & -Y\\
        -X & -Y & 1
    \end{bmatrix}\bar{\bx}_1
\end{equation}
The values of $X$ and $Y$ may be determined from the remaining two points (e.g., $\bar{\bx}_2$ and $\bar{\bx}_3$).

To solve for $X$ and $Y$, recognize that the quadratic form $\bar{\bx}_i^T \bC \bar{\bx}_i = 0$ from Eq.~\eqref{eq:ConicLocusConstraint} may be rewritten as
\begin{subequations}
\label{eq:ConstraintsP2P3}
\begin{equation}
\label{eq:ConstraintsPoint2}
0 = \bx_2^T\bC\bx_2=\left(U_1X+U_2Y+U_3\right)\left(U_4X+U_5Y+U_6\right) = \left( \bu^T_{123} \bxi \right) \left( \bu^T_{456} \bxi \right) 
\end{equation} 
\begin{equation}
\label{eq:ConstraintsPoint3}
0 = \bx_3^T\bC\bx_3=\left(V_1X+V_2Y+V_3\right)\left(V_4X+V_5Y+V_6\right)
= \left( \bv^T_{123} \bxi \right) \left( \bv^T_{456} \bxi \right)
\end{equation}
\end{subequations}
where $\bxi = [X;Y;1]$, $\bu_{ijk} = [U_i;U_j;U_k]$, $\bv_{ijk} = [V_i;V_j;V_k]$, and 
\begin{equation}\label{eq:U15}
   U_{1,4}=(r_2\pm r_1)( r_2x_1\mp r_1x_2)\qquad  U_{2,5}=(r_2\pm r_1)(r_2y_1\mp r_1y_2)
\end{equation}
\begin{equation}\label{eq:U36}
    U_3=U_6=r_1^2-r_2^2
\end{equation}
\begin{equation}\label{eq:V15}
     V_{1,4}=(r_3\pm r_1)( r_3x_1\mp r_1x_3)\qquad V_{2,5}=(r_3\pm r_1)(r_3y_1\mp r_1y_3 )
\end{equation}
\begin{equation}\label{eq:V36}
    V_3=V_6=r_1^2-r_3^2
\end{equation}
In the expressions for $U_{1,4}$, $U_{2,5}$, $V_{1,4}$, and $V_{2,5}$, the upper sign (in the $\pm$ or $\mp$ on the right-hand side) is always associated with the first of the two subscripts. 

Consider the point $\bxi = [X,Y,1] \in \PP^2$ and that $\bell^T \bxi = 0$ when the point $\bxi$ is on the line $\bell$. Thus, $\bu_{ijk}$ and $\bv_{ijk}$ may be viewed as lines on which $\bxi$ must lie. Since one of the terms on the right-hand side of both Eq.~\eqref{eq:ConstraintsPoint2} and Eq.~\eqref{eq:ConstraintsPoint3} must be zero, there are four possible combinations of lines. The solution for $\bxi$ is the intersection of these two lines.

Recalling that the intersection of two lines in homogeneous coordinates is given by their cross-product \cite{Hartley:2003}, the four possible solutions for $\bxi$ are
\begin{equation}
    \label{eq:FourCrossProducts}
    \bu_{123} \times \bv_{123}
    \quad \text{or} \quad
    \bu_{456} \times \bv_{123}
    \quad \text{or} \quad
    \bu_{123} \times \bv_{456}
    \quad \text{or} \quad
    \bu_{456} \times \bv_{456}
\end{equation}
It is shown in the appendix that the solution must be the cross product $\bu_{123} \times \bv_{123}$, so there is no need to ever compute the others. Letting $\bs = [s_1;s_2;s_3] = \bu_{123} \times \bv_{123}$, then $\bxi \propto \bs$. Since $\bxi = [X;Y;1]$, the solution for $X$ and $Y$ is simply
\begin{equation}
    X = \frac{s_1}{s_3}
    \quad \text{and} \quad 
    Y=\frac{s_2}{s_3}
\end{equation}
With $X$ and $Y$ found, the value of \(Z\) may be computed by Eq.~\eqref{eq:Z}.

\subsection{Algorithm Summary}
Having developed the theory in detail, this section provides a summary of the essential steps of the algebraic solution to the Gibbs problem. Listed here are only the steps that an analyst would require to compute a solution, as most of the equations required to formally derive the algorithm are superfluous during actual implementation.

Given three position vectors \(\{\br_i\}_{i=1}^3\), compute the basis vectors of the intermediate frame $N$ as 
\begin{equation}\label{eq:w}
    \hat{\bw} =  \frac{\br_1\times \br_2}{\|\br_1\times \br_2\|}
\end{equation}
\begin{equation}\label{eq:e2}
    \be_2 = \frac{\bw\times \br_1}{\|\bw\times \br_1\|}
\end{equation}
\begin{equation}\label{eq:e1}
    \be_1 = \be_2 \times \bw
\end{equation}
These vectors generate the rotation matrix
\begin{equation}
    \bT_N^I = \begin{bmatrix}
        \be_1^T\\
        \be_2^T\\
        \bw^T
    \end{bmatrix}
\end{equation}
that can be used to obtain the 2-D coordinates within the orbit plane
\begin{equation}\label{eq:xiyi}
    \bar{\bx}_i=\begin{bmatrix}
        x_i\\
        y_i\\
        0
    \end{bmatrix} = \bT_N^I \br_i\qquad i=1,\,2,\,3
\end{equation}
It is now possible to find
\begin{equation}\label{eq:u123v123}
    \bu_{123} = (r_1+r_2)\begin{bmatrix}
    r_2x_1-r_1x_2\\
    r_2y_1-r_1y_2\\
    r_1-r_2
    \end{bmatrix}\qquad \bv_{123}= (r_1+r_3) \begin{bmatrix}
        r_3x_1-r_1x_3\\
        r_3y_1-r_1y_3\\
        r_1-r_3
    \end{bmatrix}
\end{equation}
that may be used to calculate
\begin{equation}
    \bs = \bu_{123}\times \bv_{123} 
\end{equation}
from which
\begin{equation}\label{eq:XYZ}
    X = \frac{s_1}{s_3}\qquad Y=\frac{s_2}{s_3}\qquad Z^2 = \frac{1}{r_1^2}\bar{\bx}^T_1 \begin{bmatrix}
        -Y^2& XY & -X\\
        XY & -X^2 & -Y\\
        -X & -Y & 1
    \end{bmatrix}\bar{\bx}_1
\end{equation}
All the parameters necessary to calculate the orbit are now available. The periapsis direction in the initial frame is given by
\begin{equation}\label{eq:p_hatXYZ}
    \hat{\bp} = \frac{1}{X^2+Y^2}\bT_I^N\begin{bmatrix}
        X\\
        Y\\
        0
    \end{bmatrix}
\end{equation}
while the parameter and the semi-major axis may be evaluated as
\begin{equation}\label{eq:pXYZ}
    p = \frac{1}{\sqrt{X^2+Y^2+Z^2}}
\end{equation}
\begin{equation}\label{eq:aXYZ}
    a = \frac{\sqrt{X^2+Y^2+Z^2}}{Z^2}
\end{equation}
\section{Numerical example}
\begin{figure}[ht!]
    \centering
    \includegraphics[width=0.7\textwidth]{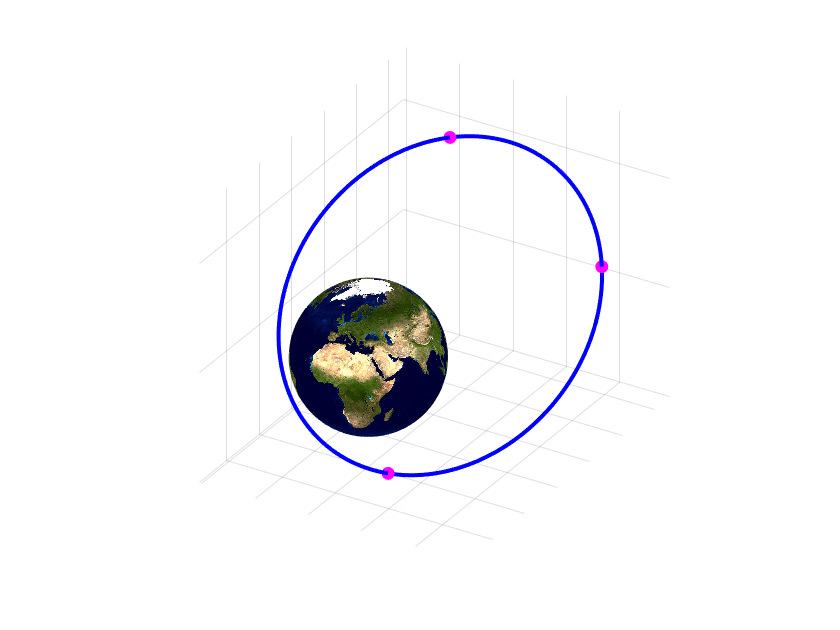}
    \caption{In blue, the orbit of the satellite. In magenta, the three position measurements used to perform IOD.}
    \label{fig:orbit}
\end{figure}

Consider the orbit shown in Fig.~\ref{fig:orbit}, described by
\[
a = 15000\,km\qquad e = 0.5\qquad i = 70\deg\qquad \Omega = 150\deg\qquad \omega =200\deg 
\]
where \(i\), \(\Omega\) and \(\omega\) are the inclination, right ascension of the ascending node and argument of the periapsis. 

Assume measurements are taken at $\{\theta_i\}_{i=1}^3 = \{70.00,\, 165.91,\, 216.49\}$ degrees, which produces the three position measurements
\[
\br_1 = \begin{bmatrix}
    1642.9 \\ 2845.6 \\ -9027.6
\end{bmatrix} \text{~km,}
\quad
\br_2 = \begin{bmatrix}
    -19201 \\ 10197 \\ 2114.2
\end{bmatrix} \text{~km,}
\quad
\br_3 = \begin{bmatrix}
    -11678 \\ 547.76 \\ 14739
\end{bmatrix} \text{~km,}
\]

Following the method of Bate, Mueller, and White \cite{BMW} for solving the Gibbs problem, it is possible to compute
\[
N = \begin{bmatrix}
    2.2536\\
    3.9034\\
    1.6405
\end{bmatrix}\times 10^{12}\,km^3 \qquad D = \begin{bmatrix}
    2.0032\\
    3.4697\\
    1.4582
\end{bmatrix}\times 10^{8}\,km^2\qquad S = \begin{bmatrix}
    -0.2889\\
    0.9579\\
   -1.8824
\end{bmatrix}\times 10^{8}\,km^2
\]
that may be used to determine the perifocal basis vectors
\[
\hat{\bp} = \begin{bmatrix}
     0.8723\\
   -0.3685\\
   -0.3214
\end{bmatrix}\qquad \hat{\bq} = \begin{bmatrix}
    -0.1355\\
    0.4493\\
   -0.8830
\end{bmatrix}\qquad \hat{\bw} = \begin{bmatrix}
    0.4698\\
    0.8138\\
    0.3420
\end{bmatrix}
\]
the semi-latus rectum and the eccentricity
\[
p = 11250\,km\qquad e = 0.5
\]
This concludes Gibbs solution.

Moving to the approach described in this note, the measurements may be used to determine the basis vectors using Eq.~\eqref{eq:w}-\eqref{eq:e1}:
\[
\hat{\bw} = \begin{bmatrix}
    0.4698\\
    0.8138\\
    0.3420
\end{bmatrix}\qquad
\be_2 = \begin{bmatrix}
    -0.8660\\
    0.5000\\
    0.0000
\end{bmatrix}\qquad\be_1 = 
\begin{bmatrix}
    0.1710\\
    0.2962\\
   -0.9397
\end{bmatrix}
\]
from which one may use Eq.~\eqref{eq:xiyi} to calculate 
\[
\begin{bmatrix}
    x_1\\
    y_1
\end{bmatrix}=\begin{bmatrix}
    9607.1\\
    0
\end{bmatrix}\qquad \begin{bmatrix}
    x_2\\
    y_2 
\end{bmatrix}=\begin{bmatrix}
    -2249.9\\
    21727
\end{bmatrix}\qquad \begin{bmatrix}
    x_3\\
    y_3
\end{bmatrix} = \begin{bmatrix}
    -15685\\
    10387
\end{bmatrix}
\]
Substituting these 2-D coordinates into Eq.~\eqref{eq:u123v123} provides 
\[
\bu_{123}=\begin{bmatrix}
    7.2795\times 10^{12}\\
    -6.5646\times 10^{12}\\
    -3.8482\times 10^{8}
\end{bmatrix}\qquad \bv_{123} = \begin{bmatrix}
     9.419\times 10^{12}\\
     -2.8361\times 10^{12}\\
     -2.6162\times 10^{8}
\end{bmatrix}
\]
At this point, \(X\), \(Y\) and \(Z\) may be calculated with Eq.~\eqref{eq:XYZ}
\[ X = 1.5201 \times 10^{-5},  
\quad 
Y = -4.1764 \times 10^{-5}, 
\quad \text{and} \quad
Z = 7.6980 \times 10^{-5} 
\]
Once \(X\), \(Y\) and \(Z\) are known, the periapsis direction may be evaluated with Eq.~\eqref{eq:p_hatXYZ}
\[
\hat{\bp} = \begin{bmatrix}
 0.8723\\
   -0.3685\\
   -0.3214
\end{bmatrix}
\]
as well as the remaining orbital elements (using Eq.~\eqref{eq:aXYZ}, Eq.~\eqref{eq:ecc} and Eq.~\eqref{eq:pXYZ}):
\[
a = 15000\,km\qquad e = 0.5\qquad p = 11250\,km
\]
As expected, these agree exactly with the original problem and with the vector-based Gibbs solution.

\section{Conclusions}
This work presents an algebraic geometry interpretation of orbit determination (OD) from three position vectors, commonly known as ``Gibbs problem.'' Central to this algorithm is a new solution to the three-point ellipse fitting problem with a constrained focus location (as opposed to the five-point fitting of an ellipse in general position). Constraining a focus to be at the origin requires that the orbit's conic envelope is a circle, which permits a simple OD solution. Indeed, after a bit of algebraic manipulation, the parameters of the orbit's conic locus (or conic envelope) may be found by the intersection of two lines. Transformation from conic locus to orbital elements is straightforward, as both are complete representations of a Keplerian orbit.

Despite the length of the derivation, only a few equations are needed when performing OD in practice. A short algorithm summary and numerical example illustrates that the computations required by the (classical) vector analysis approach and the new algebraic approach are comparable. The numerical example also illustrates that the same answer is produced by both solution methods.

\appendix
\section{Identifying the correct conic matrix }\label{app:proof}
This appendix demonstrates why $\bu_{123} \times \bv_{123}$ is the correct cross-product from Eq.~\eqref{eq:FourCrossProducts} to construct the conic matrix $\bC$. The other three options always represent a physically impossible geometry.

To begin, recall that the radius of an orbit with semi-latus rectum $p$ and eccentricity $e$ is given by \cite{BMW}
\begin{equation}\label{eq:polar}
    r = \frac{p}{1+e cos(\theta)}
\end{equation}
where \(\theta\) is the true anomaly. It is not true, however, that any conic can be completely parameterized by Eq.~\eqref{eq:polar}. In the case of a hyperbola, Eq.~\eqref{eq:polar} only describes one branch---with the other branch being described by 
\begin{equation}\label{eq:other_branch}
    r = \frac{p}{1-ecos(\theta)}
\end{equation}

The orbit parameters \(p\), \(e\) and \(\theta\) depend on \(X\), \(Y\) and \(Z\). It is possible to show that the \(X\), \(Y\) and \(Z\) values computed from selecting the incorrect cross-products in Eq.~\eqref{eq:FourCrossProducts} lead to the observed points lying on both branches of a hyperbolic conic. Since a orbit can only travel along a single branch of a hyperbola, those solutions may be discarded.

This can be proved starting from the expressions of \(p\) and \(e\) provided in section~\ref{sec:an_algebraic_view}, and from the definition of true anomaly.

Since the true anomaly is the angle between $\hat{\bp}$ and $\br$, it follows that $\cos \theta_i = \hat{\bp}^T (\br/r)$. In the orbit plane, this is 
\begin{equation}
    \cos \theta_i = \frac{1}{r_i} \hat{\bp}^T\begin{bmatrix}
        x_i\\
        y_i\\
        0
    \end{bmatrix}
\end{equation}
Substituting this for $\cos \theta_i$ in Eq.~\eqref{eq:polar}, along with  Eq.~\eqref{eq:SemiLatusRectum} and Eq.~\eqref{eq:ecc}, the right-hand side of Eq.~\eqref{eq:polar} may be rewritten as
\begin{equation}
    r_i = \frac{1}{\sqrt{X^2+Y^2+Z^2}+\frac{Xx_i+Yy_i}{r_i}}
\end{equation}
which is the same as
\begin{equation}\label{eq:den_factor}
    r_i \sqrt{X^2+Y^2+Z^2} + Xx_i+Yy_i=1
\end{equation}
Using the expressions of \(X\), \(Y\) and \(Z\) given by the choice of $\bu_{123} \times \bv_{123}$, a few manipulations lead to
\begin{equation}
  {X^2+Y^2+Z^2}=  \left(\frac{U_1 V_2 - U_2 V_1  - x_i (U_2 V_3 - U_3 V_2) + y_i (U_1 V_3  - U_3 V_1)}{r_1 (U_1 V_2 - U_2  V_1)}\right)^2 
  \end{equation}
With the expressions for \(\{U_i\}_{i=1}^3\) and \(\{V_i\}_{i=1}^3\) in Eq.~\eqref{eq:U15}-\eqref{eq:V36}, the previous equation can be compactly written as 
\begin{equation}
    X^2+Y^2+Z^2 = \left(\frac{K_2}{K_1}\right)^2
\end{equation}
where 
\begin{equation}
    K_1 =  det\begin{bmatrix}
       x_1 & y_1 & r_1\\
       x_2 & y_2 & r_2\\
       x_3 & y_3 & r_3
    \end{bmatrix}\qquad
    K_2 = det \begin{bmatrix}
        x_1 & y_1 & 1\\
        x_2 & y_2 & 1\\
        x_3 & y_3 & 1
        \end{bmatrix}
\end{equation}
so that
\begin{equation}\label{eq:BoverA}
    r_i\sqrt{X^2+Y^2+Z^2} = r_i\frac{|K_2|}{|K_1|}
\end{equation}
Similarly, it is possible to express 
\begin{equation}\label{eq:xX+yY}
    x_iX+y_iY = \frac{x_i(U_2 V_3 - U_3 V_2) - y_i(U_1 V_3 - U_3 V_1)}{U_1V_2 - U_2 V_1} = \frac{K_1-K_2r_i}{K_1}
\end{equation}
Substitution of Eq.~\eqref{eq:BoverA} and Eq.~\eqref{eq:xX+yY} into Eq.~\eqref{eq:den_factor} yields the helpful relation
\begin{equation}
    \frac{K_1|K_1| + K_1|K_2|r_i - K_2|K_1|r_i}{K_1|K_1|}=1 \implies K_1|K_2| = K_2|K_1|
\end{equation}
It becomes clear, then, that Eq.~\eqref{eq:polar} is satisfied for all \(i\) when \(K_1\) and \(K_2\) have the same sign.

Repeating the process with the expressions of \(X\), \(Y\) and \(Z\) associated with \(\bu_{456}\times \bv_{123}\), Eq.~\eqref{eq:den_factor} may be rewritten for \(i=1\), \(i=2\) and \(i=3\). After manipulations, the two equations become:
\begin{equation} \label{eq:K3r1}
    \frac{K_3|K_3| - K_2|K_3|r_1 + K_3|K_2|r_1}{K_3|K_3|}=1 \implies K_2|K_3| = K_3|K_2|
\end{equation}
\begin{equation} \label{eq:K3r2}
    \frac{K_3|K_3|+K_2|K_3|r_2 + K_3|K_2|r_2}{K_3|K_3|}=1 \implies K_2|K_3| = -K_3 |K_2|
\end{equation}
\begin{equation} \label{eq:K3r3}
    \frac{K_3|K_3| - K_2|K_3|r_3 + K_3|K_2|r_3}{K_3|K_3|}=1 \implies K_2|K_3| = K_3|K_2|
\end{equation}
where 
\begin{equation}
    K_3 = det\begin{bmatrix}
    x_1 & y_1 & r_1\\
    x_2 & y_2 & -r_2\\
    x_3 & y_3 & r_3
    \end{bmatrix}
\end{equation}

The conditions on the right-hand side of Eq.~\eqref{eq:K3r1} and Eq.~\eqref{eq:K3r2} cannot be true at the same time. This contradiction occurs because the procedure assumed all the points were on the branch of the hyperbola described by Eq.~\eqref{eq:polar}, whereas the choice of \(\bu_{456}\times \bv_{123}\) actually places the three points on both branches of the hyperbola (points $\br_1$ and $\br_3$ are one branch, while point $\br_2$ is on the other). 

\begin{figure}[ht!]
    \centering
    \includegraphics[width=0.5\textwidth]{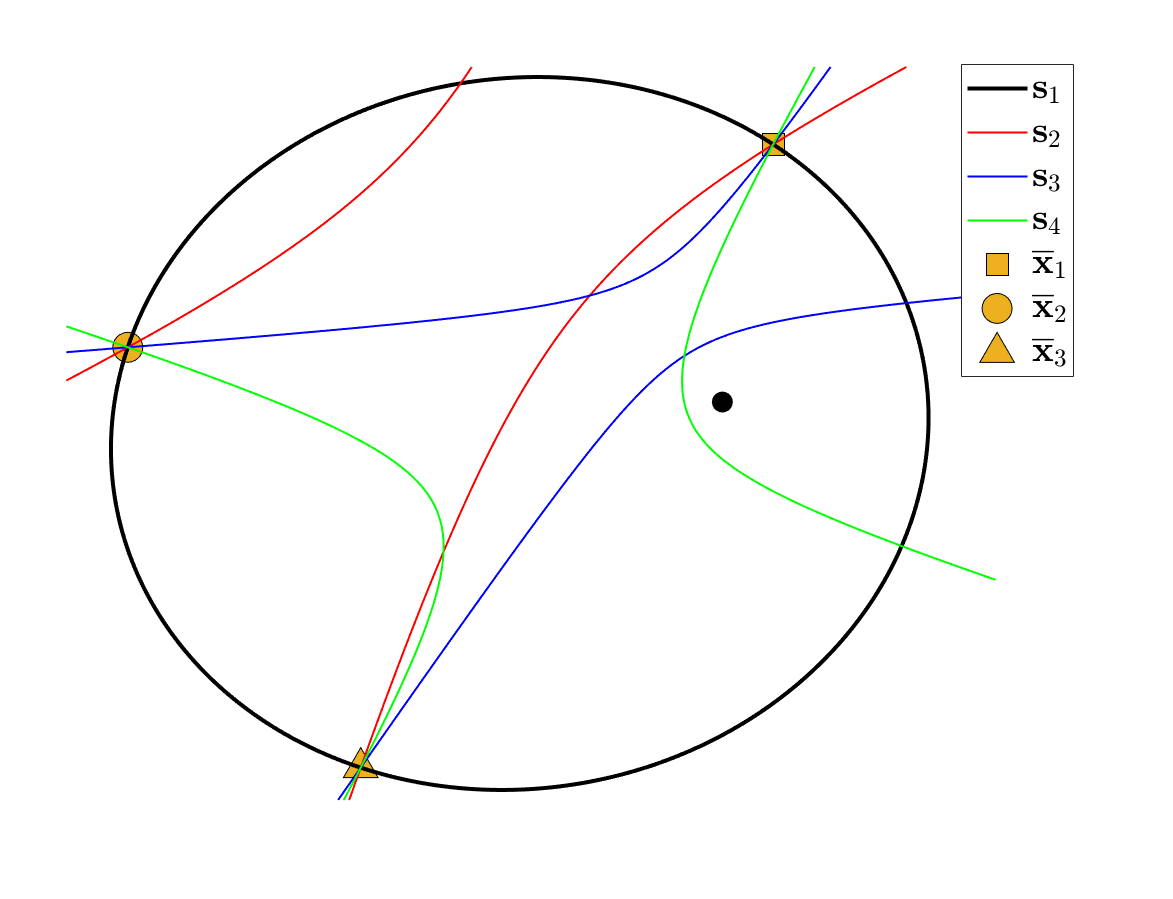}
    \caption{The conics corresponding to the solutions \(\bs_{2}\), \(\bs_3\) and \(\bs_4\) always have one of the points lying on a branch different from the other two.}
    \label{fig:example1}
\end{figure}

The same procedure may be used to show that the conics corresponding to \(\bu_{123}\times \bv_{456}\) and \(\bu_{456}\times \bv_{456}\) are not acceptable as well, and \(\bu_{123}\times \bv_{123}\) remains the only solution that does not impose contradictory conditions. An example of this behavior is shown in Fig.~\ref{fig:example1}.

\section*{Acknowledgment}
The authors thank Liam Smego for reviewing a draft of this manuscript and providing helpful feedback.

\bibliography{sample}

\end{document}